# A new method for the study of boundary value problems for linear and quasi-linear systems of ODE

**Prof. Konyaev Yu. A.**

**Peoples Friendship University of Russia**

St. Mikluho-Maklaya, B. 13, Moscow, 117198, Russia

**Department of Mathematics**

**(Moscow)**

*The method is proposed for the study of many-point boundary value and spectral problems for systems of ODE, by reducing them to special equivalent integral equations, and allows us [in contrast with the known method [1]] to consider boundary and initial value problems. Here we don't use the mechanism of Green's function whose construction is quite nontrivial, especially in case of many-point boundary value problems. The proposed algorithm [2, 3] makes it easier to study and write out the condition of unique solvability of such problems.*

The method is proposed for the study of many-point boundary value problems for systems of ordinary differential equations, by reducing them to special equivalent integral equations, and allows us [in contrast with the known method [1]] to consider boundary and initial value problems. Here we avoid the mechanism of Green's function whose construction is quite nontrivial, especially in case of many-point boundary value problems. The proposed algorithm [2, 3] makes it easier to study and write out the condition of unique solvability of such problems.

**Theorem 1.** The many point boundary value problem in $\mathbb{R}^n$:

$$\dot{x} = A(t)x + f(t) \; ; \quad \sum_1^m F_j x(t_j) = \alpha \quad (2 \le m \le n) \tag{1}$$

$$(x, f \in \mathbb{R}^n \; ; \; A(t) \; ; \; f(t) \in C[0,1] \; ; \; 0 = t_1 < t_2 < \ldots < t_m = 1)$$

where the $F_j$'s are constant matrices under the condition $det F_j \neq 0$ $(F = \sum_1^m F_j \Phi(t_j)$ ; $\Phi(t)$ is any fundamental matrix) have a unique solution, and its representation is given by the form:

$$x(t) = \Phi(t)C + \sum_1^m \Phi_k(t) \int_{t_k}^{t} \Phi^{-1}(s) \qquad (2)$$

where

$$C = F^{-1}\left(\alpha - \sum_{k=1}^{m} F_j \sum_{k=1}^{m} \Phi_k(t_j) \int_{t_k}^{t_j} \Phi^{-1} f(s) ds\right) ;$$

$$\dot{\Phi}_k = A(t)\Phi_k, \quad (k = \overline{1,m}) \quad ; \quad \sum_{k=1}^{m} \Phi_k(t) = \Phi(t)$$

**Proof.** It follows from the fact that the term $\Phi(t)C$ is the general solution of the homogenous $(f(t) \equiv 0)$ system (1), and the sum $\sum_1^m \Phi_k(t) \int_{t_k}^{t} \Phi^{-1}(s) f(s) ds$ represents the particular solution of the non-homogenous system (1), which can be verified directly.

**Remarks.** 1. The integral representation (2) coincides with the known Cauchy formula for the solution of initial value problems when $m = 1$, and if $(2 \leq m \leq n)$ then the integral representation is a solution for different variants of boundary value problems.

2. The given assertion is constructive under the fundamental matrix $\Phi(t)$. If for some "near" matrix $B(t)$ holds that $\Psi(t)(\dot{\Psi}(t) = B(t)\Psi)$ and the difference $(\|A(t) - B(t)\|)$ is small in norm, then under the condition that $det F_0 \neq 0$ $(F_0 = \sum_1^m F_j \Psi(t_j))$ problem (1) reduces to an integral equation of the form

$$x(t) = \Psi(t)C + \sum_1^m \Psi_k(t) \int_{t_k}^{t} \Psi^{-1}(s)[(A(s) - B(s))x + f(s)] ds \equiv Lx$$

and under the contraction operator $Lx$ problem (1) is uniquely solvable.

3. If $f(x,t)$ is nonlinear and nonhomogeneous, then we get the integral equation (2).

The proposed algorithm for analyzing many-point boundary value problems ( Theorem 1), and taking into account the splitting method [2, 3] allows us to study more complex classes of singularly perturbed boundary-value problems of the form:

$$\varepsilon \dot{x} = A(t,\varepsilon)x + f(t,\varepsilon) \ ; \quad (x, f \in \mathbb{R}^n) \qquad (3)$$

$$(\sum_{1}^{n} F_j x(t_j, \varepsilon) = a \quad t \in [0,1] \ ; \ 0 < t_1 < t_2 < \ldots < t_n = 1 \ ; \ n \geq 2, 3)$$

[ where $F_j$ denotes certain constant matrices $(j = \overline{1,n})$ and the function series $f(t,\varepsilon) = \sum_{1}^{\infty} f_k(t)\varepsilon^k$ and the matrix series $A(t,\varepsilon) = \sum_{1}^{\infty} A_k(t)\varepsilon^k$ are from sufficiently smooth vector functions $f_k(t)$ and matrix functions $A_k(t)$ respectively and are convergent absolutely and uniformly under a certain norm on the interval [0,1] for sufficiently small ε (0 <ε<<1) ].

The term singularity in this case means the solution of the limit (ε = 0) of ($3_0$) may not satisfy the boundary conditions.

In some regard, the structure of the spectrum $\{\lambda_{0j}(t)\}_1^n$ of the matrix $A_0(t)$ with the position of points $t_k$ and certain other conditions [2, 3] of the singularly perturbed problem (3) is uniquely solvable and its solution will be contained in the neighborhood of $t_k$ $(k = \overline{1,n})$ in the boundary layers of exponential type [2, 3], provided that the spectrum $\{\lambda_{0j}(t)\}_1^n$ of the matrix $A_0(t)$ satisfies the inequalities:

$$\lambda_{0j}(t) \neq \lambda_{0k}(t) \ ; \ (j \neq k \ ; \ j,k = \overline{1,n} \ ; \ t \in [0,1]) \ ; \lambda_{0j}(t) \neq 0 \qquad (4)$$

$$Re\lambda_{01}(t) \leq 0 \ ; \ Re\lambda_{0n}(t) \geq 0 \ ; \quad t \in [0,1]$$

$$Re\lambda_{0k}(t) \geq 0 \ ; \ t \in [0, t_k] \ , \ Re\lambda_{0k}(t) \leq 0 \ ; \ t \in [t_k, 1] , \ (k = \overline{2, n-1})$$

$$Re \int_{t_k}^{t} \lambda_{0k}(s)ds < 0 \ ; \ t \in [0,1]\backslash t_k \quad (k = \overline{1,n}).$$

In applying condition (4), the quasi-regular asymptotic solution of the singularly perturbed problem (3) will be written with any required accuracy in closed analytic form and is suitable for the qualitative and numerical analysis.

It should be noted that for the analysis of singularly perturbed boundary value problems other asymptotic methods [4- 5] are not applicable or they are not effective.